\documentclass[12pt]{amsart}
\usepackage{amssymb,amscd,latexsym,eufrak,amsmath,graphicx,epsfig}
\usepackage{amsthm}

\theoremstyle{plain}
\newtheorem{Theorem}{Theorem}
\newtheorem{Corollary}[Theorem]{Corollary}

\newtheorem{Lemma}[Theorem]{Lemma}
\newtheorem{Proposition}[Theorem]{Proposition}

\theoremstyle{definition}
\newtheorem{Definition}[Theorem]{Definition}
\newtheorem{Notation}[Theorem]{Notation}
\newtheorem{Remark}[Theorem]{Remark}

\newtheorem{Example}[Theorem]{Example}

\newcommand{\ceiling}[1]{{\lceil #1 \rceil}}
\newcommand{\floor}[1]{{\lfloor #1 \rfloor}}
\newcommand{\fractional}[1]{{\{ #1 \}}}
\def\ring(#1){{{\mathcal R}(#1)}}
\def\locus(#1){{{\mathcal L}(#1)}}
\newcommand{\ga}{{\mathfrak a}}

\newcommand{\tensor}{\otimes}

\newcommand{\CC}{{\mathbb C}}
\newcommand{\RR}{{\mathbb R}}

\newcommand{\OO}{{\mathcal O}}

\newcommand{\multr}[1]{{\mathcal J}(r \cdot #1)}

\newcommand{\mult}[1]{{\mathcal J}(#1)}

\begin{document}

\title{Multiplier Ideals of Sufficiently General Polynomials}

\author{Jason Howald}

\thanks{Jason Howald, Department of Mathematics, University of Michigan,
Ann Arbor, MI, 48109, \email{jahowald@umich.edu}}



\thanks{Thanks to Rob Lazarsfeld for many helpful discussions, and for his patience.}

\begin{abstract}
It is well known that the multiplier ideal $\multr{\ga}$ of
an ideal $\ga$ determines in a straightforward way 
the multiplier ideal $\multr{f}$
of a sufficiently general element $f$ of $\ga$.
We give an explicit condition on a polynomial 
$f \in \CC[x_1,\ldots,x_n]$ which guarantees that it is a 
sufficiently general element of the most natural associated monomial
ideal, the ideal generated by its terms.  This allows us to directly  
calculate the multiplier ideal $\multr{f}$ (for all $r$) 
of ``most'' polynomials $f$.
\end{abstract}

\maketitle

\section{Introduction}

The multiplier ideal of a divisor \(D\) or ideal sheaf \(\ga\)
tells us a great deal about the singularities of \(D\), or of a
general section of \(\ga\).  It is also a main ingredient of the
Nadel vanishing theorem, which of course has found numerous
important applications. The multiplier ideal is a powerful tool
because it provides ``invariants with vanishing.''  In this paper
we will compute the multiplier ideals of a large class of divisors
in ${\mathbb A}^n$, those defined by ``nondegenerate''
polynomials. \footnote{see \cite{Arnold} for earlier uses of
nondegeneracy} This nondegeneracy condition is fairly weak: It is
satisfied, for example, by functions of the form $c_1x_1^{a_1} +
c_2x_2^{a_2} + \ldots + c_nx_n^{a_n}$.

\begin{Notation}
We write $\multr{\ga}$ for the multiplier ideal of an ideal $\ga$,
with constant $r$.  We write $\multr{D}$ for the multiplier ideal
of a divisor $D$.  For $f \in \CC[x_1,\ldots,x_n]$, we abbreviate
$\multr{Div(f)}$ by ``$\multr{f}$,'' although ``${\mathcal
J}(f^r)$'' would be more suggestive.
\end{Notation}

We will not develop the basic properties of multiplier ideals.  We
refer the reader instead to \cite{PAG} for a comprehensive
algebraic development.

The algebraic definition of the multiplier ideal refers to an
embedded resolution of \(D\) or \(\ga\), making it difficult to
calculate directly in specific examples.  In a previous paper
\cite{HowaldMonomials}, we calculated the multiplier ideal of an
arbitrary monomial ideal \(\ga\) in a polynomial ring
$\CC[x_1,\ldots,x_n]$. It happens that if \(\ga \subset
\CC[x_1,\ldots,x_n]\) is an ideal and \(f \in \ga\) is a generic
section\footnote{meaning a generic \(\CC\)-linear combination of
generators \(\{g_1,\ldots,g_k\}\) of \(\ga\).} of \(\ga\), then
the multiplier ideals of \(Div(f)\) and $\ga$ are essentially the
same:
\begin{Theorem}
Let $\ga \subset \CC[x_1,\ldots,x_n]$  be a monomial ideal, and
let $f \in \ga$ be a generic section.  If $r<1$ then
\[\multr{(f)}=\multr{\ga}.\]
\end{Theorem}

See \cite[Proposition 9.2.29]{PAG} for a proof.  If $r \geq 1$,
the relationship is only slightly more complicated.

This begs the question:  When is a given function a ``generic
enough'' element of some related monomial ideal?  We give a
sufficient condition which guarantees that a polynomial \(f \in
\CC[x_1,\ldots,x_n]\) is a sufficiently generic element of the
(monomial) ideal $\tau(f)$ generated by its terms.  This allows us
to calculate the multiplier ideal of \(f\).

Roughly, this condition requires that for every face of the Newton
polyhedron of $f$, the function obtained by summing just those terms of $f$ which
lie on that face should have nonvanishing derivative on the
torus $T^n \subset {\mathbb A}^n$.  (For each face $\sigma$ of the
Newton polyhedron, we call this ``$\sigma$-nondegeneracy.'')  We
prove that this guarantees equality of $\multr{f}$ and
$\multr{\tau(f)}$ in four steps:

\begin{enumerate}
\item  There is a toric log resolution $\mu:X' \rightarrow X$, which
principalizes $\tau(f)$.

\item  By toric geometry, $X'$ is locally isomorphic to ${\mathbb
A}^n$, so we may regard $\mu^*(f) = f'$ as a polynomial.  We prove
that it inherits nondegeneracy from $f$.

\item  Since $f'$ is nondegenerate and $\tau(f')$ is principal,
$Div(f')$ is a normal crossing divisor, and $\mu$ is a
simultaneous toric log resolution for both $\tau(f)$ and $Div(f)$.

\item  Direct comparison of the terms $r\mu^*(Div(f))$ and
$r\mu^{-1}(\tau(f))$ appearing in the definition of the multiplier
ideal shows them to be the same.


\end{enumerate}

\section{Main Theorem and Examples}

The nondegeneracy condition which we will use describes the
relationship of a polynomial \(f\) with the faces of its Newton
polyhedron. We must first introduce some definitions and notation
for dealing with such structures.

\begin{Definition}
Let \(X=\CC^n\).  In questions of membership in monomial ideals, the
coefficient of a monomial \(m \in \CC[x_1,\ldots,x_n]\) is
irrelevant.  Discarding this information, we will refer to the
lattice \(L_X \cong {\mathbb Z}^n\) of rational monomials on
\(X\), whose nonnegative orthant contains the monomials.  Let
\(\ga \subset \OO_X\) be a monomial ideal.  We may identify
\(\ga\) with a subset of \(L_X\).  The {\it Newton polyhedron}
\(P(\ga)\) of \(\ga\) is the convex hull of this set, an unbounded
closed subset of \(L_X \tensor \RR\).  By a {\it face} \(\sigma\)
of the Newton polyhedron, we mean a face of arbitrary codimension
(even the codimension-$0$ face) . Because the polyhedron is
unbounded, faces need not be compact.
\end{Definition}

\begin{Definition}
Let \(f \in \CC[x_1,\ldots,x_n]\) be a polynomial.  The {\it term
ideal} \(\tau(f)\) of \(f\) is the monomial ideal generated by the
terms of \(f\).  We write \(P(f)\) for \(P(\tau(f))\), though this
is not the usual definition of the Newton polyhedron of a
polynomial.
\end{Definition}

We can now introduce the main nondegeneracy condition which will
guarantee that \(f\) and \(\tau(f)\) have similar multiplier
ideals.

\begin{Definition}
Let \(f \in \CC[x_1,\ldots,x_n]\), and let \(\sigma\) be a face of
\(P(f)\). We write \(f_\sigma\) for the polynomial composed of the
terms of \(f\) which lie in \(\sigma\).  Note that \(f_\sigma\) remembers
the coefficient of each term.
If the \(1\)-form
\(df_\sigma\) is nonvanishing on the torus \(T_n=(\CC \setminus
0)^n\), then we say that \(f\) is {\it nondegenerate for
\(\sigma\)}. If \(f\) is nondegenerate for every face of its
Newton polyhedron, then we call \(f\) simply {\it nondegenerate}.
If \(f\) is nondegenerate for just its compact faces, then we say
{\it \(f\) has nondegenerate principal part}.
\end{Definition}

This last condition follows \cite[6.2.2]{Arnold}, in which it is
used (from the analytical point of view) in the calculation of the
log canonical threshold of \(f\). Notice that if \(\sigma\) is the
origin in \(L_X\), then \(f_\sigma\) is a nonzero constant.  In
this special case, we call \(f\) nondegenerate for \(\sigma\)
despite the vanishing of \(df_\sigma\).  At the other extreme,
\(\sigma\) may be the entire Newton polyhedron, in which case
nondegeneracy for \(\sigma\) implies that \(Div(f|_{T^n})\) is
nonsingular.  (If this condition seems too strong, remember that
we will allow ourselves only toric resolutions, and these are
isomorphisms on $T^n$.)

We will soon see that nondegeneracy of \(f\) implies convenient
global properties for \(Div(f)\), whereas nondegeneracy of the
principal part tells us that \(Div(f)\) is ``nice'' only near the
origin.

\begin{Definition}
Let \(X = \CC^n\), and let \(\ga\) be an ideal sheaf on \(X\).  By
a {\it toric log resolution} of \(\ga\), we mean a log resolution
\(\mu:X' \rightarrow X\) which is obtained by a sequence of
blowings-up along orbits of the torus action.  Equivalently, it is
a log resolution for which $X'$ and $\mu$ belong to the toric
category.  We define a toric log resolution of a divisor
similarly.
\end{Definition}
Locally on \(X'\), such a map is determined by a monomial ring map
\(\CC[x_1,\ldots,x_n] \rightarrow \CC[y_1,\ldots,y_n]\).  A toric
log resolution of \(\ga\) can also be described in terms of
refinements of the dual polyhedron of \(P(\ga)\).  See
\cite{HowaldMonomials} for more on this perspective.

\begin{Proposition}
Every monomial ideal has a toric log resolution.
\end{Proposition}

See \cite{HowaldThesis} for a sketch of the proof, or
\cite{Fulton} for the basic theory of toric varieties.

\begin{Example}

The functions \(f= y^2-y(x-1)^2\), \(g = (xy-1)^9\), and
\(h=(x+y)^2-(x-y)^5\) on \({\mathbb C}^2\) are all degenerate.
Although \(df\) is nonzero on the torus, there is a proper face
\(\sigma \subset P(f)\) for which \(f_\sigma= y(x-1)^2\).  Clearly
\(df_\sigma=0\) at the point \((1,0)\). \(Div(f)\) has a
singularity at \((1,0)\), but toric resolutions up can only blow
up the origin so the singularity cannot be corrected by such maps.
Nevertheless, \(f\) has nondegenerate principal part.

For every proper face \(\sigma\) of \(P(g)\), \(g_\sigma\) is
nondegenerate.  But \(dg\) itself vanishes at points in the torus,
so \(g\) is degenerate.  Although \(Div(g)\) has a toric log
resolution, the high multiplicity of \(g\) along the proper
transform of \(Div(g)\) complicates the calculation of the
multiplier ideal.

For properly chosen \(\sigma\), \(h_\sigma\) is a perfect square
\((x+y)^2\), and is clearly degenerate.  Despite being singular
only at the origin, \(h\) has no toric log resolution.  The
multiplier ideal of \(Div(h)\) is not monomial, even for some
\(r<1\).  This function has degenerate principal part.
\end{Example}

\begin{Example} \label{Diagonal divisor is nondegenerate}
    Let \(f = c_1x_1^{a_1}+\ldots+c_nx_n^{a_n}\), with \(c_i \neq 0\) and \(a_i > 0\).
    Then \(f\) is nondegenerate, as the reader may check.
    Furthermore, if the simplex spanned by \(\{a_ie_i\}_{i \in \{1\ldots n\}}\)
    (the unique maximal compact face of \(P(f)\)) contains no
    lattice points other than its vertices, then any polynomial
    with the same Newton polyhedron as \(f\) has nondegenerate principal
    part.
\end{Example}

\begin{Notation}
We write \(\multr{\ga}\) or \(\multr{D}\) for the multiplier ideal
of an ideal or divisor, with rational coefficient \(r\).  We also
write \(\floor{r}\) and \(\fractional{r}\) for the round-down and
fractional part of a rational number \(r\).  We will write
\(\floor{D}\) for the round-down of a \({\mathbb Q}\)-divisor
\(D\), and \(\ceiling{D}\) for its round-up.
\end{Notation}

The multiplier ideal of a monomial ideal can be calculated
combinatorially:

\begin{Theorem} \label{MonomialMult}
    Let \( {\mathfrak a} \subset {{\mathcal O}_{{\mathbb A}^{n}}}\) be a monomial
    ideal.  Then \(\multr{{\mathfrak
    a}} \) is a monomial ideal, and contains exactly the following
    monomials:
    \[ \{m : m x_1x_2 \ldots x_n \in Interior(rP(\ga))\}.\]
\end{Theorem}

\begin{proof}
    See \cite{HowaldThesis} or \cite{HowaldMonomials}.
\end{proof}

For small values of \(r\), we can hope for the equality
\(\multr{Div(f)}=\multr{\tau(f)}\), provided that \(f\) is a
sufficiently generic element of \(\tau(f)\).  This requires a
Bertini argument which shows that a toric log resolution \(\mu\)
of \(\tau(f)\) is in fact a log resolution of a general section of
\(\tau(f)\).  The divisors \(\mu^*(Div(f))\) and \(\mu^{-1}(\ga)\)
differ only by the proper transform of \(Div(f)\), which
disappears on rounding if \(r<1\).

So a claim that \(f\) is sufficiently general primarily means that
any toric log resolution of \(\tau(f)\) also resolves \(Div(f)\).
Our result is in essence an effective substitute for the use of
Bertini's theorem to guarantee this.  Basic facts about multiplier
ideals will take care of the rest of the calculation. Here then is
our main theorem.

\begin{Theorem}[Main Theorem] \label{MAIN}
Let \(f \in \CC[x_1,\ldots,x_n] \).  Let \(\mu:X' \rightarrow
\CC^n\) be a toric log resolution of \(\tau(f)\).  If \(f\) is
nondegenerate, then \(\mu\) also log-resolves \(Div(f)\).  If
\(f\) has nondegenerate principal part, then there is a Zariski
neighborhood \(U\) of the origin over which \(\mu\) is a log
resolution of \(f|_U\).
\end{Theorem}

We will prove this theorem in the next section.

\begin{Corollary}\label{nondeg mult r<1}
Let \(f \in \CC[x_1,\ldots,x_n]\) be nondegenerate, and let
\(r<1\). Then
\[\multr{Div(f)} = \multr{\tau(f)}.\]
\end{Corollary}
\begin{proof}[Proof (\ref{MAIN} \(\Rightarrow\) \ref{nondeg mult r<1})]
    Since we have a toric log resolution \(\mu:X' \rightarrow X=\CC^n\)
    of both \(f\) and \(\tau(f)\),
    we may use it to calculate each of these multiplier ideals.
    Let \(K_{X'/X}\) be the relative canonical bundle, and
    consider the divisors
    \[F_1=K_{X'/X} - \floor{r\mu^*Div(f)}
    \hspace{0.2in} \text{ and } \hspace{0.2in}
    F_2=K_{X'/X} - \floor{r\mu^{-1}(\tau(f))},\]
    which when pushed forward via \(\mu_*\) calculate the two
    multiplier ideals.  We claim that these divisors are equal.
    Because the resolution is toric, each exceptional divisor
    \(E\) of \(\mu\) corresponds to a weighted blowup.  The
    order along \(E\) of \(\mu^*Div(f)\) is the minimal weight of
    the terms of \(f\), and the order along \(E\) of
    \(\mu^{-1}(\tau(f))\)
    is the minimal weight of monomials in \(\tau(f)\).
    These are obviously the same, so \(F_1\) and \(F_2\) have the
    same order along each exceptional divisor.  Clearly
    they also agree on the proper transforms of the coordinate
    axes.

    We have tested equality on all divisors not
    intersecting the torus \(T_n\).
    No other divisors can appear in the support of \(F_2\)
    because \(\tau(f)\) is monomial.  We must similarly rule out components
    of \(F_1\) which intersect the torus.  Suppose \(D'\) were one such.
    Since \(\mu\) is an isomorphism on the torus, \(\mu(D')\) would be a
    divisor \(D\) and
    \(ord_{D'}(\mu^*(f)) = ord_D(f) = 1\) because
    \(df|_{T_n}\) does not vanish.  To calculate \(ord_{D'}F_1\),
    we would multiply this number by \(r\) and round down, giving
    \(0\) because \(r<1\).
\end{proof}

\begin{Corollary} \label{NondegenerateMult}
Let \(f \in \CC[x_1,\ldots,x_n]\) be nondegenerate.  Then
\[\mult{r \cdot Div(f)} = (f^\floor{r})\mult{\fractional{r} \tau(f)}.\]
\end{Corollary}
\begin{proof}
    The integer divisor \(\floor{r}Div(f)\) pulls outside the
    multiplier ideal (see \cite{HowaldThesis}), so we have
    \[\mult{r \cdot Div(f)} = (f^\floor{r})\mult{\fractional{r} Div(f)}.\]
    Corollary \ref{nondeg mult r<1} finishes the proof.
\end{proof}

\begin{Example}
    As in Example \ref{Diagonal divisor is nondegenerate}, let \(f =
    c_1x_1^{a_1}+\ldots+c_nx_n^{a_n}\), with \(c_i \neq 0\) and \(a_i > 0\).
    Let \(v = (\frac{1}{a_1},\ldots,\frac{1}{a_n})\), so that
    \(m \in P(f)\) if and only if \(v \cdot m \geq 1\).  (Here we regard \(m\) as
    a vector in \(L_X\) and take its dot product with \(v\).)
    Theorem \ref{MonomialMult} tells us that
    \(\multr{\tau(f)}\) is generated by \(\{m:
    v \cdot (m x_1x_2 \ldots x_n) > r\}\).  By Corollary \ref{NondegenerateMult},
    \[\multr{Div(f)} \text{ is generated by } \{f^\floor{r} m:
    v \cdot (m x_1x_2 \ldots x_n) > \fractional{r}\}\]
\end{Example}

\section{Proof of Main Theorem}

We will require a few more definitions and notations before we can
comfortably prove the main theorem.  Let \(X=\CC^n\).

\begin{Definition}
Let \(\ga \subset \CC[x_1,\ldots,x_n]\) be a monomial ideal, and
let \(\sigma \subset P(\ga)\) be a face of its Newton polyhedron.
Let \(v_\sigma\) be any linear functional on \(L_X\) which cuts
out \(\sigma\) in the sense that \(v_\sigma(P(\ga)) = [b,\infty)\)
and \(v_\sigma^{-1}(b) \cap P(\ga) =\sigma\). We may regard \(v\)
as a hyperplane distinguishing \(\sigma \subset P(\ga)\).  Let
\(\ring(\sigma)\) be the ring generated by all monomials \(m\)
with \(v_\sigma(m)=0\).  Because \(v_\sigma\) is a nonnegative
vector (when written in the standard dual basis),
\(\ring(\sigma)\) is of the form \(\CC[S]\) for some \(S \subseteq
\{x_1,\ldots,x_n\}\). Let \(\locus(\sigma)=
Spec(\ring(\sigma))
\subset \CC^n\). Notice that \(\locus(\sigma)\) is just a
coordinate linear subspace of \(\CC^n\).  By the {\it relative
interior} of \(\locus(\sigma)\), we will mean the points of
\(\locus(\sigma)\) not contained in some smaller coordinate linear
subspace.
\end{Definition}

The reader may check that although \(v_\sigma\) is not completely
determined by \(\sigma\), \(\ring(\sigma)\) and \(\locus(\sigma)\)
are.  Readers familiar
with valuations will notice that \(v_\sigma\) is a valuation on
the ring \(\CC[x_1,\ldots,x_n]\), and \(\locus(\sigma)\) is its center.

\begin{Example} \label{LocusExample}
    Let \(\ga= (x^3,y^3,z^3) \subset \CC[x,y,z]\).  The
    triangular face \(\sigma\) of \(P(\ga)\) is determined by
    \(v_\sigma=(1/3,1/3,1/3)\).  For this face and in fact for all
    of its subfaces, \(v_\sigma\) is strictly positive, so
    \(\ring(\sigma)=\CC\) and
    \(\locus(\sigma)\) is the origin.
    On the other hand, if we take \(\sigma\) to be the
    intersection of \(P(\ga)\) with the \(x\)-axis in the lattice of monomials,
    then \(\ring(\sigma)=\CC[x]\), and \(\locus(\sigma)\) is the
    \(x\)-axis in \(\CC^3\).
\end{Example}

\begin{Proposition} \label{OriginImpliesCompact}
    If, as in Example \ref{LocusExample}, \(\locus(\sigma)\) is
    the origin, then \(\sigma\) is compact.
\end{Proposition}
\begin{proof}
If \(\locus(\sigma)\) is the origin, then \(\ring(\sigma)=\CC\),
so \(v_\sigma\) is a strictly positive dual vector.  From the
above definition, \(v_\sigma^{-1}(b) \cap P(\ga) =\sigma\).  Since
\(v_\sigma^{-1}(b)\) has compact intersection with the nonnegative
orthant in \(L_X\), and since \(P(\ga)\) is closed, \(\sigma\)
itself must be compact.
\end{proof}

\begin{Definition} \label{PushforwardOfFace}
    Let \(X=Spec \ \CC[x_1,\ldots,x_n]\), let \(X'=Spec \ \CC[y_1,\ldots,y_n]\),
    and let \(\mu:X' \rightarrow X\) be a
    monomial map.  If \(\ga \subset \OO_X\) is a monomial ideal,
    then so is \(\ga'=_{def}\mu^{-1}\ga\).  Let \(\sigma'\) be a face of
    \(P(\ga')\).  We write \(\mu(\sigma')\) for
    \(\{m \in L_X : \mu^*(m) \in \sigma'\}\).  It is
    the preimage of \(\sigma'\) under \(\mu^*\).
\end{Definition}

\begin{Proposition} \label{PushforwardIsAFace}
    For any face \(\sigma'\) of \(P(\ga')\),
    \(\mu(\sigma')\) is a face of \(P(\ga)\).
\end{Proposition}
\begin{proof}
    Whereas \(\sigma'\) is the intersection of
    \(P(\ga')\) with the level set of some linear
    functional \(v_{\sigma'}\), \(\sigma =_{def}\mu(\sigma')\) is the
    intersection of \(P(\ga)\) with the level set of \(v =_{def} v_{\sigma'}
    \circ \mu^*\).  This \(v\) is a linear functional because \(\mu^*\)
    acts linearly on the lattice of monomials.  The reader may
    verify that \(v(\sigma)\) is the smallest value of \(v\) on
    \(P(\ga)\).  This proves that \(\sigma\) is a face.
\end{proof}

\begin{Proposition} \label{CommuteMuLocus}
    Let \(X\), \(X'\), \(\ga\), \(\ga'\), and \(\mu:X' \rightarrow X\) be as
    in Definition \ref{PushforwardOfFace}.
    For any face \(\sigma'\) of \(\ga'\),
    \(\mu(\locus(\sigma'))=\locus(\mu(\sigma'))\).
\end{Proposition}
\begin{proof}
    Since everything is invariant under the torus action, each set
    is at least torus invariant, and so is defined by a monomial ideal in
    \(\CC[x_1,\ldots,x_n]\).
    Choose \(v'=v_{\sigma'}\) cutting out \(\sigma'\).  As in
    Proposition \ref{PushforwardIsAFace}, \(\sigma\) is cut out
    by \(v = v' \circ \mu^*\).  Unraveling the
    definitions involved, we find that
    the ideal of the left hand side is generated by
    \(\{m \in L_X: v'(\mu^*(m))>0\}\).  The ideal of the
    right hand side is generated by \(\{m \in L_X: (v' \circ
    \mu^*)(m)>0\}\).  These are obviously the same
    set.
\end{proof}

We now have the tools in hand to prove the main theorem, which we
restate here:

\begin{Theorem}[Main Theorem]
Let \(f \in \CC[x_1,\ldots,x_n] \).  Let \(\mu:X' \rightarrow
\CC^n\) be a toric log resolution of \(\tau(f)\).  If \(f\) is
nondegenerate, then \(\mu\) also log-resolves \(Div(f)\).  If
\(f\) has nondegenerate principal part, then there is a Zariski
neighborhood \(U\) of the origin over which \(\mu\) is a log
resolution of \(f|_U\).
\end{Theorem}

We begin with a special case:

\begin{Lemma} \label{PrincipalSubcase}
    Let \(f \in \CC[x_1,\ldots,x_n] \), and assume that
    \(\tau(f)\) is principal.  Let \(\sigma\) be a face of
    \(P(f)\).  If \(f\) is nondegenerate for \(\sigma\), then
    \(Div(f)\) is a normal crossing divisor at points \(p\) in
    the relative interior of \(\locus(\sigma)\).
\end{Lemma}
\begin{proof}

    Let \(m\) be the (unique) monomial generator of \(\tau(f)\).
    By reordering the
    coordinates, we may assume \(\sigma=\{mx_1^{a_1}\cdot \ldots
    \cdot x_k^{a_k}:a_i \geq 0\}\) for some \(k\), so
    \(\locus(\sigma)=Spec(\CC[x_1,\ldots,x_k])={\mathbb V}(x_{k+1},\ldots,x_n)\).
    We may write \(g=f/m\), where
    \(m\) is a monomial and \(g\) has nonzero constant term.  Let
    us also write \(g_\sigma=f_\sigma/m\).  Notice that \(g=g_\sigma\) on
    \(\locus(\sigma)\), and \(g_\sigma \in \CC[x_1,\ldots,x_k]\).
    (See the figure below.)

\begin{figure}[h]
\centerline{\epsfig{file=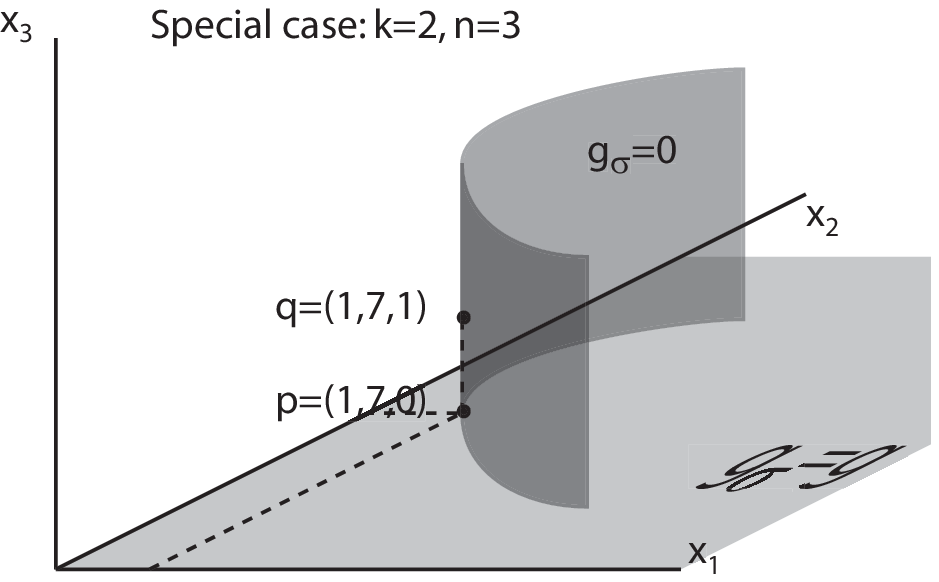}}
\end{figure}

    Let \(p\) be in the relative interior of \(\locus(\sigma)\).
    If \(g(p) \neq 0\), the normal crossing conclusion is trivial, so we
    assume \(g(p)=0\).
    We must show that the relevant \(1\)-forms
    \(\{dx_{k+1},\ldots,dx_n,dg\}\) corresponding to the components of
    \(Div(f)\) which meet at \(p\) are independent at \(p\).
    For this, it suffices to show that \(d(g|_{\locus(\sigma)})\)
    is nonzero at \(p\).  Let \(q = p+(0,\ldots,0,1,\ldots,1)\) (k
    zeroes, n-k ones).  Notice that \(g_\sigma(q)=g_\sigma(p)=g(p)=0\).
    Since \(q\) is in the torus, we have at \(q\),
    \[0 \neq df_\sigma = d(mg_\sigma) = m dg_\sigma+g_\sigma dm = m dg_\sigma.\]
    So \(0 \neq dg_\sigma\) at \(q\).
    Since \(g_\sigma \in \CC[x_1,\ldots,x_k]\),
    \[ 0 \neq dg_\sigma |_q = dg_\sigma |_p = d(g_\sigma|_\locus(\sigma))
    |_p= d(g|_\locus(\sigma))|_p,\]
    as desired.
\end{proof}

\begin{Remark}
    Notice that every point \(p \in \CC^n\) is in the relative
    interior of some coordinate linear space \(L\), and such a
    space can always be written \(L=\locus(\sigma)\) for some
    \(\sigma\).
\end{Remark}

\begin{proof}[Proof of Theorem \ref{MAIN}]
    Let \(f'=\mu^*(f)\).
    We must show that \(Div(f')\) has normal crossing
    support.  Since this property is local on \(X'\), and since
    \(\mu\) is a {\it toric} log resolution of \(\tau(f)\), we may take
    \(X'= Spec\ \CC[y_1,\ldots,y_n]\), with \(\mu:X' \rightarrow
    X\) a monomial map.  Note that \(\tau(f')\) is principal,
    because \(\mu\) resolves \(\tau(f)\).

First we discuss the sense in which \(f'\) inherits nondegeneracy from \(f\).
Let  \(\sigma'\) be a face of \(P(f')\), and let $\sigma = \mu(\sigma')$
be the associated face of \(P(f)\).
Clearly \(\mu^*f_\sigma = f'_{\sigma'}\).  Because $\mu$ is an isomorphism
on the torus, $df'_{\sigma'}|_T = df_\sigma|_T$.  This shows that $f'$ is
$\sigma'$-nondegenerate as long as $f$ is $\sigma$-nondegenerate, where
$\sigma = \mu(\sigma')$.


    If we assume that \(f\) is nondegenerate for every face then
    we can conclude the same for \(f'\), so \(Div(f')\)
    is a normal crossing divisor (by  Lemma \ref{PrincipalSubcase})
    and we are done.

    If however we assume only that \(f\) is nondegenerate for its compact
    faces, then we know only that \(f'\) is nondegenerate for each
    \(\sigma'\) with \(\mu(\sigma')\) compact.  In particular (by
    Proposition \ref{OriginImpliesCompact}), \(f'\) is
    nondegenerate for each \(\sigma'\) with
    \(\locus(\mu(\sigma'))\) equal to the origin.  Equivalently,
    (by Proposition \ref{CommuteMuLocus}), for each \(\sigma'\) with
    \(\mu(\locus(\sigma'))\) the origin.  We may conclude (again by
    Lemma \ref{PrincipalSubcase}) that
    \(f'\) has normal crossing support near points \(p' \in X'\)
    which map to the origin.  Since the locus of point \(p'\)
    which fail the normal crossing condition is a closed subset
    of \(X'\), and since \(\mu\) is proper, there is a
    neighborhood of the origin above which \(\mu\) is a log
    resolution of \(Div(f)\).  Therefore if $f$ has nondegenerate
    principal part, then $\mu$ is indeed a log resolution of $f$,
    over a neighborhood of the origin.
\end{proof}

\bibliographystyle{apalike}
\bibliography{citations}

\end{document}